\documentclass[10pt,bezier]{article}
\usepackage{amsmath}
\usepackage{amsfonts,amsthm,amssymb}
\usepackage{amsfonts}
\usepackage{graphics}
\usepackage{amscd}
\usepackage{amsmath,amsfonts,amssymb,amscd}
\usepackage{indentfirst,graphicx,epsfig}
\usepackage{graphicx,psfrag}
\input{epsf}

\begin{document}
\title{\bf Rainbow connection of bridgeless outerplanar graphs with small diameters\thanks{
The project supported partially by CNNSF (No.11271204, No.11301371), Central South University of Forestry and Technology Project (No.(104$\vert$0292)0292). Deng XC
supported by Tianjin Normal University Project (No.52XB1206).
}}
\author{Xingchao Deng$^a$, He Song$^b$, Guifu Su$^{c}$, Weihua Yang$^{d}$  \\
 \footnotesize  $^a$College of Mathematical Science, Tianjin Normal University\\
   \footnotesize  Tianjin City, 300387, P. R. China\\
        \footnotesize  E-mail: dengyuqiu1980@126.com\\
 \footnotesize $^b$School of Mathematical Sciences, LPMC, Nankai University\\
    \footnotesize Tianjin City, 300071, P. R. China\\
       \footnotesize E-mail: songhe@mail.nankai.edu.cn\\
\footnotesize $^c$School of Science, Beijing University of Chemical Technology\\
  \footnotesize  Beijing, 100029, P. R. China\\
 \footnotesize $^d$Department of Mathematics, Taiyuan  University of Technology\\
  \footnotesize  Taiyuan, 030024, P. R. China\\}
     %\footnotesize  E-mail: ywh222@163.com\\}
 %\footnotesize $^e$College of Sciences, Central South University of Forestry and Technology\\
%  \footnotesize  Changsha, 410004, P. R. China }
\date{}
\maketitle
\begin{abstract}
In this paper, we investigate rainbow connection number $rc(G)$ of bridgeless outerplanar graphs $G$ with diameter 2 or 3. We proved the following results:
If $G$ has diameter $2,$ then $rc(G)=3$ for fan graphs $F_{n}$ with $n\geq 7$ or $C_5,$ otherwise $rc(G)=2;$ if $G$ has diameter $3,$ then $rc(G)\leq 4$ and the bound is sharp.\\

{\bf MSC2010 subject classifications:} 05C55, 05C12.\\


\end{abstract}

\section{Introduction}

\vspace{3mm}

Throughout this paper, we consider finite undirected graphs without loops
and multiple edges. For any graph $G$, whose vertex and edge sets are denoted by $V(G)$ and
 $E(G)$ respectively. We define the order of $G$ by $n=n(G)=\mid V(G)\mid$ and
the size by $m=m(G)=\mid E(G)\mid$. For terminology and notations not defined here we refer to $\cite{WT}.$ The {\emph diameter}, denoted by $diam(G)$, of a graph $G$ is the maximal distance between any two vertices; the {\emph radius}, denoted by $rad(G)$, of a graph $G$ is the minimum eccentricity. A subset $D$ of $V(G)$ is called a {\emph dominating set} if every vertex
of $G-D$ has at least one neighbor in $D,$ and further $D$ is called a connected dominating
set if the {dominating set} $D$ induces a connected subgraph of $G.$
Let $X, Y\subseteq V(G),$ we say $X$ dominates $Y$ if every vertex of $Y$ is
adjacent to at least one vertex of $X.$ The {\emph join} of two graphs $G$ and $H$ is a graph which is obtained by linking each vertex of $G$ to all vertices in $H$.
We use $K_{n}$ and $P_{n}$ to denote the complete graph and path of order $n$ respectively, and $K_{s,t}$ to denote the complete bipartite graph with each partite having $s$ and $t$ vertices, respectively.

\vspace{1mm}

\vspace{1mm}

Recall an {\emph outerplanar graph} is a graph that can be embedded in a plane without crossings in such a way that all the vertices lie in the boundary of the unbounded embedded face.

\vspace{1mm}

The following is a characterization of outerplanar graphs due to \cite{GF}.

\vspace{1mm}

\noindent\textbf{Theorem A} \cite{GF}. A graph is outerplanar if and only if it has no $K_4$ or $K_{2,3}$ minor.
\vspace{1mm}

We call a path in an {\emph edge-colored connected} graph $G$ is a rainbow path if edges of it are colored different. In 2008, Chartrand et al.\ $\cite{CJMZ}$ introduced the concept of rainbow connection of graphs. An edge-colored graph $G$
is said to be {\emph rainbow connected} if there is a rainbow path between every pair of vertices.
The {\emph rainbow connection number} of $G$, denoted by $rc(G)$, is the minimum number of colors needed to make $G$ rainbow connected. The rainbow connection number of any complete graph is 1 and that of a tree equals to its size.

\vspace{1mm}

Note that cut-edges must have different colors when $G$ is rainbow connected, thus stars have arbitrarily large rainbow connection number while having diameter $2.$ Therefore, it is interesting to seek an upper bound of rainbow connection number in terms of its diameter for 2-edge-connected graphs. Chandran et al.\ $\cite{CDRV}$ studied the relation between rainbow connection numbers and connected dominating set, and showed:
$\mathbf{(i)}$ For any bridgeless chordal graph $G,$ $rc(G)\leq 3 rad(G)$;
$\mathbf{(ii)}$ for any unit interval graph $G$ with minimum degree at least $2,$ $rc(G)=diam(G).$
The upper bound is sharp. Moreover they proved that $rc(G)\leq rad(G)(rad(G)+2)$ if $G$ is 2-edge-connected. Li and co-authors $\cite{LLL}$ proved that $rc(G)\leq 5$ for any 2-edge-connected graph $G$ with diameter 2.
Dong $\cite{JD}$ showed the upper bound 5 is tight. Li et al.\ $\cite{LLS}$ also presented that $rc(G)\leq 9$ for all 2-edge-connected graph $G$ with diameter 3.

\vspace{1mm}

Recently Huang et al.\ $\cite{XH1}$investigated upper bounds on the rainbow connection number of  bridgeless outerplanar graphs. In particular, they proved that:
If $diam(G)=2,$ then $rc(G)\leq 3$ and the bound is tight; if $diam(G)=3,$ then $3\leq rc(G) \leq 6.$

\vspace{1mm}

Let $\mathcal{G}_{n}^{D}=\{G|G~\text{is a bridgeless outerplanar graphs with order}~n~\text{and}~diam(G)=D\}.$ In this paper, we give the rainbow connection number for graphs of $\mathcal{G}_{n}^{2}$ and obtained a tight upper bound 4 for the rainbow connection number of graphs in
$\mathcal{G}_{n}^{3}.$
\vspace{3mm}

\section{Rainbow connection number of $G\in \mathcal{G}_{n}^{2}$}

\vspace{3mm}

 We begin with the following exact formulas of rainbow connection number for special graphs.

\vspace{1mm}

\noindent\textbf{Theorem 2.1} (\cite{CJMZ}). For cycle $C_n$, we have
\begin{equation*}
rc(C_n)=\left\{
\begin{array}{ll}
 \frac{n}{2} & {\rm if~}n\ {\rm is\ even,}\\
\lceil\frac{n}{2}\rceil & {\rm if}~n~{\rm is\ odd.}
\end{array}
\right.
\end{equation*}

\noindent\textbf{Theorem 2.2} (\cite{XD}). For fan $F_n=P_n\vee K_1,$ where $\vee$ is join operation of two graphs, we have
  \begin{equation*}
rc(F_n)=\left\{
\begin{array}{ll}
1 & {\rm if}~~n=2,\\
2 & {\rm if}~~3\leq n\leq 6,\\
3 & {\rm if}~~n\geq 7.
\end{array}
\right.
\end{equation*}

\vspace{1mm}

We now state the main result of this section.

\vspace{1mm}

\noindent{\bf Theorem 2.3}. Let $G$ be a bridgeless outerplanar graph with order $n$ and diameter $2$. Then \begin{equation*}
rc(G)=\left\{
\begin{array}{ll}
3 & {\rm if}~~G=F_{n}~(n\geq 7)~ {\rm or~ C_5,}\\
~~\\
2 & {\rm otherwise}.
\end{array}
\right.
\end{equation*}

\vspace{1mm}

\noindent{\bf Proof.} Assume $G$ is a  bridgeless outerplanar graph with diameter $2,$ then each of induced cycle has length at most $5.$

\vspace{1mm}

For simplicity, we distinguish the following three cases:\\

\noindent{\bf Case $1.$} The longest induced cycle in $G$ is $C_5.$   In the case, it is easily verified that $G$ is  $C_5$, then $rc(G)=3.$

\vspace{1mm}

\noindent{\bf Case $2.$} The longest induced cycle in $G$ is a $C_4.$

\vspace{1mm}

Since any 4-cycle is an induced cycle with diameter 2, without loss of generality we assume that this 4-cycle is $C_{4}=a-b-c-d$. Obviously, it is exactly a bridgeless outerplanar graph with diameter $2$. Note that $C_4$  is the graph with the smallest size in the case. In the following we can construct all graphs satisfy the condition in Theorem 2.3 by adding new vertices to $C_{4}$ orderly.

{\bf Step $1.$} Adding a new vertex $e$ to $C_{4}=a-b-c-d$, bearing in mind that $G$ is a  bridgeless outerplanar graph, we know that $e$ is adjacent to at least 2 vertices among $a,b,c,d.$

\vspace{1mm}

If $e$ is adjacent to at least $3$ vertices of $a,b,c,d,$ then $G$ contains a $K_{2,3}$ minor contradicting Theorem A.
If $e$ is connected to exactly $2$ vertices, the resulting graph is exactly the graph $(1)$ and $(2)$ respectively depicted in Fig.1. Note that there exists a $K_{2,3}$ minor in $(1)$ which is also impossible. Hence, $G$ is just isomorphic to the graph $(2)$ which has rainbow connection number 2.

\vspace{1mm}

{\bf Step $2.$} Adding another new vertex $f$ to $(2)$ in Fig.1. Then $f$ is adjacent to at least 2 vertices among $a,b,c,d$ and $e.$  We will show that $f$ is adjacent to $a,b$ or $c,d$.

\vspace{1mm}

In fact, if $f$ is adjacent to 3 vertices of $a,b,c,d,e,$ then $G$ contains a $K_{2,3}$ minor, a contradiction to Theorem A.
If $f$ is connected to exactly $2$ non-adjacent vertices, again a contradiction because the existence of a $K_{2,3}$ minor. Hence, $f$ is only connected
to $2$ adjacent vertices among these five vertices. We also claim that $f$ is not adjacent to $b$ and $d$, otherwise there exists a $K_{2,3}$ minor. If $f$ is adjacent to $a,c$, then the resulting graph has diameter $3$ which contradicts the diameter being 2.  Therefore, $f$ can be connected to $a,b$ depicted in Fig.1. and similarly  can be connected to $c,d.$
\vspace{1mm}

According to the outer-planarity and diameter 2, we cannot add additional vertices to graph $G$. Easily, we can obtain $rc(G)=2.$

\vspace{1mm}

%\begin{center}
%\scalebox{0.6}{\includegraphics{diameter-2longestinducedcycle4.eps}
%$\mbox{Fig.1. diameter-2~longest ~induced~ cycle 4}$}
%\end{center}
%\begin{center}
%\scalebox{0.6}{\includegraphics{diameter-2longestinducedcycle4.eps}}
%\end{center}
%\begin{figure}[ht]

\begin{figure}[ht]
\begin{center}
\includegraphics[width=10cm,totalheight=25mm]{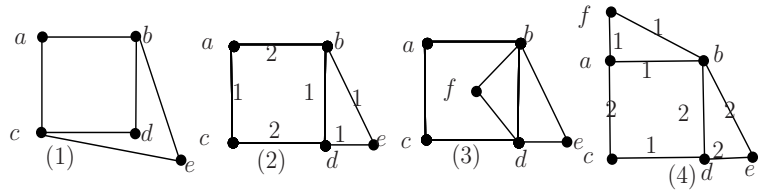}
\end{center}
\caption{$\mathrm{Diameter-2~longest ~induced~C_4}$}\label{fig1}
 \end{figure}

%\begin{center}
%\includegraphics[width=8cm]{fig1diameter2.eps}
%\end{center}
%\caption{Diameter-2~longest ~induced~4-cycle}\label{fig1}
%\end{figure}

\noindent{\bf Case 3.} The longest induced cycle in $G$ is $C_3.$

\vspace{1mm}

In this case, $G$ is a maximal outerplanar graph(MOP).
The authors $\cite{LB}$ presented that an MOP can be recursively constructed as follows: $\mathbf{(a)}$ The graph $K_3$ is an MOP graph.
$\mathbf{(b)}$ For an MOP graph $H_1$ embedded in the plane with vertices lines on the exterior face $F_1,$ let $H_2$ denoted the graph obtained by joining a new vertex to the $2$ vertices of an edge on $F_1.$ Then $H_2$ is an MOP graph. $\mathbf{(c)}$ Any MOP graph can be constructed by repeatedly applying $\mathbf{(a)}$ and $\mathbf{(b)}.$

\vspace{1mm}

In the following, we will construct outerplanar graphs under the condition of Case 3.

\vspace{1mm}

{\bf Step $1.$} We begin with $H_{1}=K_3$ depicted in Fig.2.

\vspace{1mm}

{\bf Step $2.$} Bearing in mind the symmetry of $K_3,$ we add a new vertex $d$ which is adjacent to any 2 vertices among $a,b$ and $c$. Without loss of generality, $d$ is adjacent to vertex $c$ and $b.$ The rainbow connection number of $K_{3}+d$ is 2.
%\begin{figure}[ht]
%\begin{center}
%\includegraphics[width=4cm]{fig2.eps}
%\end{center}
%\caption{$\mathrm{H_1=K_3}$}\label{fig2}
%\end{figure}
%\begin{center}
%\scalebox{0.6}{\includegraphics{H1K3.eps}}
%\end{center}

%\begin{center}
%\scalebox{0.6}{\includegraphics{fig3.eps}}
%\end{center}

\begin{figure}[ht]
\begin{center}
\includegraphics[width=10cm,totalheight=25mm]{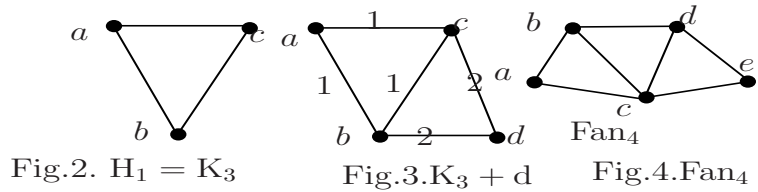}
\end{center}
\caption{$\mathrm{Three~ figures~ of~ Case3}$}\label{fig1}
\end{figure}

%\begin{figure}[ht]
%\begin{center}
%\includegraphics[width=4cm]{Three Figures of case3.eps}
%\end{center}
%\caption{$\mathrm{Three~ figures ~of ~Case3}$}\label{fig3}
%\end{figure}

%\begin{center}
%\scalebox{0.6}{\includegraphics{H1K3.eps}}
%\end{center}

{\bf Step $3.$} Since only two vertices $a$ and $d$ in $H_{2}=K_{3}+d$ has distance 2 showen in Fig.3. We now add one another vertex, say $e$, and adjacent it to 2 vertices of $V(H_{2})$. Up to isomorphism, the resulting graph is exactly the fan $Fan_{4}$ which is an outerplanar graph. By {\bf Theorem 2.2}, we know that $rc(Fan_{4})=2.$

%\begin{center}
%\scalebox{0.6}{\includegraphics{Fan4.eps}}
%\end{center}

{\bf Step $4.$} By similar argument, we add another vertex, say $f$, to the fan $Fan_{4}$. Bearing in mind the symmetry of $Fan_{4}$, all the possible resulting graphs will be $T_{1}$ and $T_{2}$ depicted in Fig.5.

\vspace{1mm}
The endpoints of every edge on the exterior face $T_{1}$ of Fig.5 have a common vertex with distance $2$ to them. Thus we could not add new vertex to $T_{1}$ to get a big MOP having diameter $2$ with more vertex. We also give $T_{1}$ a $2$ rainbow coloring, therefore $rc(T_1)=2.$ \\

The endpoints of any edge except $ca,cf$ on the exterior face $T_{2}$ have a common vertex with distance $2$ to both of them. In the process $(b)$ of the construction of MOP,  we could only connect the new vertex to $c,a$  or $c,f.$ Both operations would produce big Fans.\\

 In the construction of outerplanar graphs with diameter $2,$  if $Fan_5$ is the $H_1$ in the definition of MOP, we could obtain $Fan_n\ (n\geq 6)$ by adding a new vertex to $H_1$, which has rainbow connection number $2.$ By the similar argument, we know that $Fan_{n}\ (n\geq 7)$ are outerplanar graphs with diameter 2 and rainbow connection number $3.$

%Note that the distance between any two vertices in Fig.5 is at most two, which shows that we could not obtain a new outplanar graph such that $T_{1}$ and $T_{2}$ as its subgraphs. By calculating, we get the rainbow connection number of $T_{i}$, $i=1,2$, equals to 2.
\begin{figure}[ht]
\begin{center}
\includegraphics[width=10cm,totalheight=25mm]{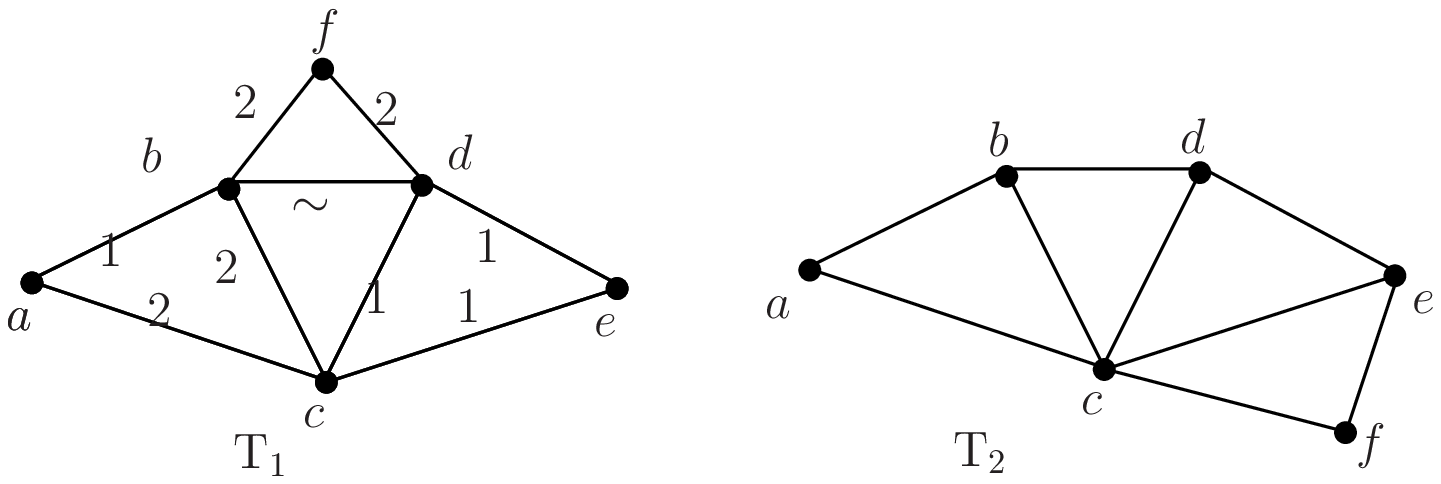}
\end{center}
\caption{$\mathrm{Fig.5. Fan_4+f}$}\label{fig1}
\end{figure}

%\begin{center}
%\includegraphics[width=4cm]{fan4+f.eps}
%\end{center}
%\caption{$\mathrm{Fig.5. Fan_4+f}$}\label{fig3}
%\end{figure}
%
%\begin{center}
%\scalebox{0.6}{\includegraphics{Fan41.eps}}
%\end{center}

This completes the proof of Theorem 2.3.
$\hfill \Box$

\vspace{1mm}

%{\it Since every edge on the exterior face $T_{1}$ of $Fig_5,$  it's endpoints have a distance $2$ vertex in $T_{1}.$  Thus we could not add new vertex to $T_{1}$ to get a big MOP with diameter $2.$\\
%
% Because any two endpoints of an edge except $ca,cf$ on the exterior face $T_{2}$  have vertex with distance $2$ to them.
% In the process $(b)$ of the construction of $MOP,$  we could only connect the new vertex to $c,a$  or $c,f.$ Thus we obtain a big $Fan.$\\
%
%
% In the construction of outerplanar graphs with diameter $2$ when $Fan_5$ is the $H_1$ in the definition of $MOP $ and added
% new vertex to $H_1$, we could obtain $Fan_n$ $(n\geq 6)$  and it has rainbow connection number $2.$ All $Fan_{n},(n\geq 7)$ have rainbow connection number
% $3.$}

\section{Rainbow connection number of $G\in \mathcal{G}_{n}^{3}$}

\vspace{1mm}

\noindent{\bf Theorem $3.1.$} Let $G$ be a bridgeless outerplanar graph with diameter $3,$ then $rc(G)\leq 4$ and the upper bound is sharp.\\

\vspace{1mm}

\noindent{\bf Proof.} We will consider  Case $1,$  Case $2,$  Case $3,$  Case $4$
and  Case $5$ to prove theorem 3.1.\\

\noindent{\bf Case $1.$} The longest induced cycle in $G$ is $C_6,$ then $G$ is $C_6$ or one of the following 5 graphs as Fig.6 shown.

\vspace{1mm}

  If $G$ is $C_6,$ its rainbow connection number is $3$ by Theorem 2.1. By the following three steps, we obtain the graphs depicted in Fig.6.

{\bf Step $1.$}  Adding a new vertex to $C_6,$ which is adjacent to at least three vertices of $C_{6}.$ It is easy to verify that the resulting graph contains $K_4$ minor, which contradicts to Theorem $A.$ Hence,the new vertex must be adjacent to 2 adjacent vertices of $C_6,$ otherwise the resulting graph will contains a $K_{2,3}$ minor. Hence, the graph obtained is exactly $(a)$ of Fig.6 under the meaning of isomorphism, where we give it a $3$ rainbow coloring. The graph has diameter $3$ obviously.  Thus $rc((a))=3.$
\vspace{1mm}

{\bf Step $2.$} Adding another new vertex to $(a)$ of Fig.6, by analogous argument as Step 1, up to isomorphism, the resulting graph is either $(b)$ or $(c)$ of Fig.6.

\vspace{1mm}

{\bf Step $3.$} We add another new vertex to the graph $(b)$ under the restriction of diameter 3 and meaning of isomorphism.  For similar reason of Step 1, the new vertex must be adjacent to adjacent
vertices of the graph $(b).$ The resulting graph must be $(e)$ of Fig.6.

 Now we add another new vertex to the graph $(c)$. The resulting graph is $(d)$ of Fig.6, by the same reason of the case: adding another new vertex to the graph $(b).$  We cannot add new vertex to graphs $(d)$ and $(e),$  otherwise the resulting graphs will have diameter $4.$

All graphs of Fig.6 have diameter $3$ and 3 rainbow coloring as depicted in it, thus their rainbow number are $3.$

%\begin{center}
%\scalebox{0.6}{\includegraphics{6-cycle+fan.eps}}
%\end{center}
%$\mbox{Fig.6.}~C_6+\mbox{Fan}$

\begin{figure}[ht]
\begin{center}
\includegraphics[width=10cm,totalheight=40mm]{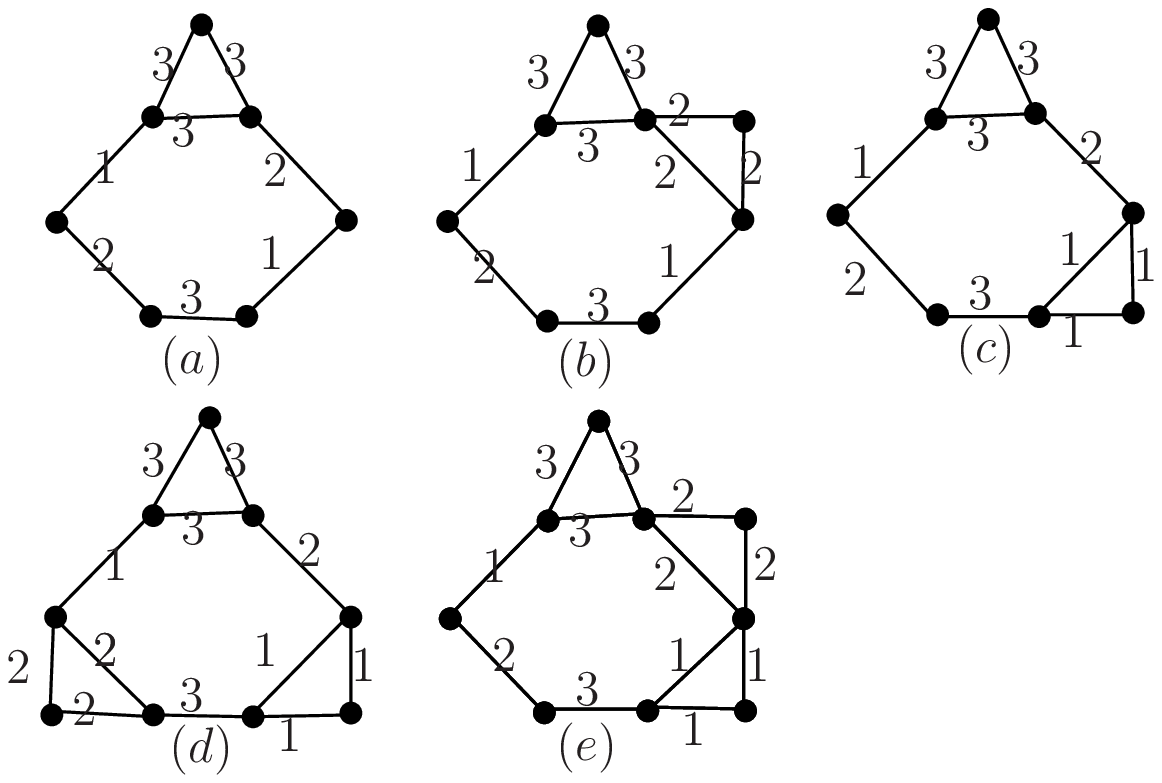}
\end{center}
\caption{$\mathrm{Fig.6.~ C_6+\mbox{Fan}}$}\label{fig1}
\end{figure}

%\begin{center}
%\includegraphics[width=4cm]{6-cycle+fan.eps}
%\end{center}
%\caption{$\mathrm{Fig.6. C_6+\mbox{Fan}}$}\label{fig3}
%\end{figure}

\noindent{\bf Case $2.$} The longest induced cycle in $G$ is $C_5$ and $G$ has an induced $C_4.$    Then $G$ are as follows.

\vspace{1mm}

Following, we add vertices one after the other to $C_5$ to obtain all graphs of Case 2.

{\bf Step $1.$} $C_5$ is preliminary of the constructing process of Case 2. When we add a new vertex with degree at least 3 to $C_5,$ the resulting graph has a $K_4$ minor, which is contrary to Theorem $A.$ Thus the new vertex must be adjacent to $2$ vertices which are adjacent in $C_5,$ otherwise the new graph will contains a $K_{2,3}$ minor. Remember that the graph has an induced $C_4.$ So the new vertex must be adjacent to $2$ vertices with distance $2$ in $C_5$; and the resulting graph has a $K_{2,3}$ minor, which is impossible. Therefore, the induced $C_4$ must be obtained by adding $2$ new vertices to the $C_5.$ The new graph must be the $(1)$ of Fig.7 by the sense of isomorphism. Moreover the graph contains only one induced $C_4$ and $C_5,$ otherwise it has diameter more than $3.$ Note $(1)$ of Fig.7 has a 3-rainbow coloring, thus $rc((1))=3.$

{\bf Step $2.$} Add new vertices to $(1)$ of Fig.7. Similar to Step 1 in Case 2, we see that the new graphs may have Fan structures with $e$ or $f$ as central vertex in the regions $A,B,C$ or $D$ of $(2)$ in Fig.7.

\begin{figure}[ht]
\begin{center}
\includegraphics[width=10cm,totalheight=30mm]{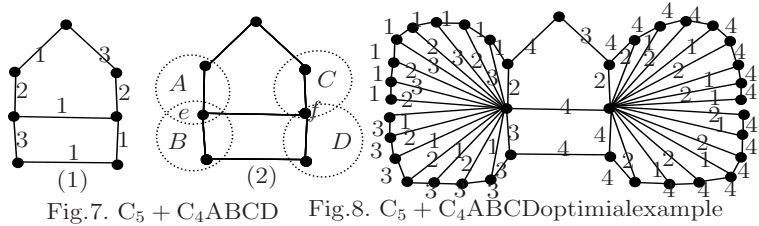}
\end{center}
\caption{$\mathrm{Two~ figures ~of ~case~ 2  }$}\label{fig1}
\end{figure}

%\begin{center}
%\scalebox{0.6}{\includegraphics{{5-cycle+4-cycleABCD.eps}}}
%\end{center}
%Two figures of case 2.eps
%\begin{center}
%\scalebox{0.6}{\includegraphics{{5-cycle+4-cycleoptimialexample.eps}}}
%\end{center}
In Fig.8, any $Fan_9$ can be replaced by another $Fan_n\ (n\geq 10\ \mbox{or}\ 1\leq n\leq 8).$
The corresponding graphs can be rainbow colored by coloring the $Fan_n$ with the same method as Fig.8 shown with $4$ colors.
All graphs have rainbow connection number at most $4,$ see Fig.8.\\

\noindent{\bf Case $3.$} The longest induced cycle in $G$ is $C_4.$

 We add vertices to $C_4$ to obtain all graphs in this case. It is easy to see that the resulting graph $G$ is one of the following cases shown as Fig.9 or Fig.10, where every Fan-structure can be replaced by any $Fan_n$ with $n\geq 1.$

When the Fan-structure is replaced by any other $Fan_n$ with $n\geq 1,$ we can color the resulting graph by the method as Fig.10 shown with $4$ colors. Since every graph in Fig.9 is an induced subgraph of Fig.10, its rainbow connection number is at most $4.$

\begin{figure}[ht]
\begin{center}
\includegraphics[width=10cm,totalheight=50mm]{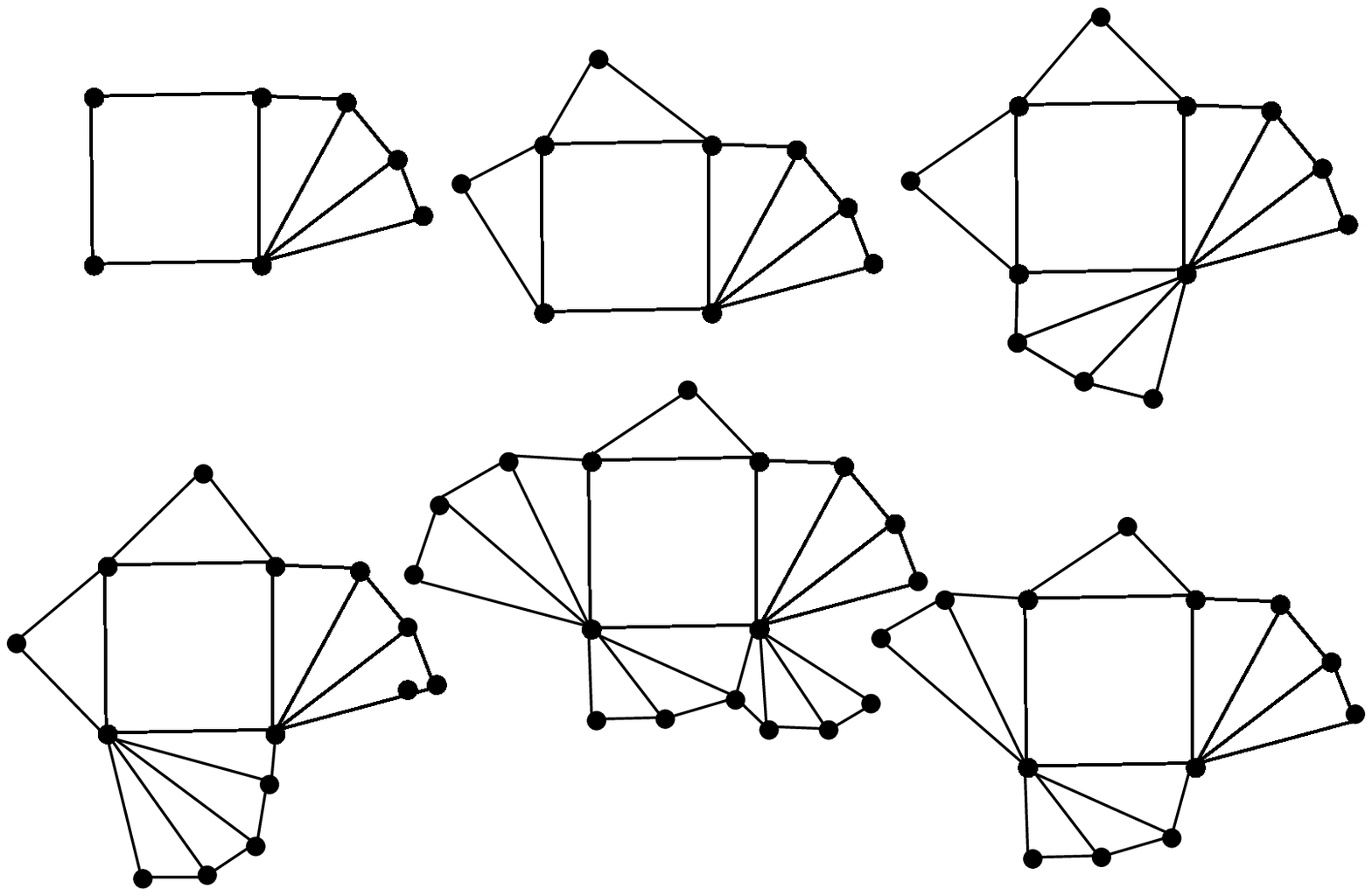}
\end{center}
\caption{$\mathrm{4-cycle+Fan  }$}\label{fig1}
\end{figure}

%\begin{center}
%\scalebox{0.6}{\includegraphics{4-cycle+Fan.eps}}
%\end{center}

%\begin{center}
%\scalebox{0.6}{\includegraphics{4-cycle+Fanoptimalexample.eps}}
%\end{center}

\begin{figure}[ht]
\begin{center}
\includegraphics[width=10cm,totalheight=37mm]{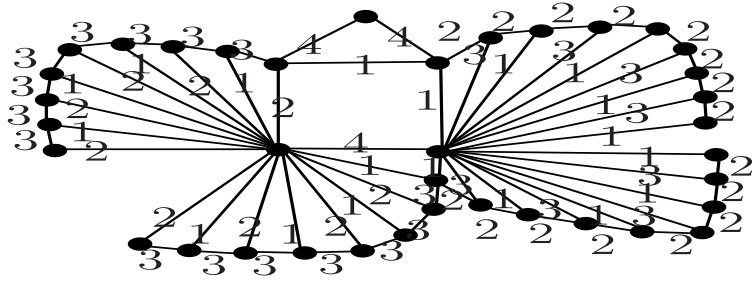}
\end{center}
\caption{$\mathrm{4-cycle+Fan ~optimal ~example }$}\label{fig1}
\end{figure}

\noindent{\bf Case $4.$} $G$ has an induced longest $C_5$ and the other induced cycles are $C_3s.$
\vspace{1mm}

 Notice the construction definition of MOPs described in Case 3 proving Theorem 2.3 and our graph has diameter $3.$ When we add new vertices to $(1)$ of Fig.11, we can obtain a big graph having Fan-structures with $d,e$ or $f$ as their central vertices. Then all graphs for Case 4 are the following ones shown as $(2)$ of Fig.11 or Fig.12.

Fig.12 gives some graphs having Fan-structure with $f$ as central vertex and at most $7$ vertices. Moreover they have rainbow coloring with at most $4 $ colors as shown in Fig.12.
The graph $(2)$ of Fig.11 gives an example with $d,e$ as the central vertices of Fan-structures ($f$ cannot be a central vertex of the Fan-structures in this setting, otherwise the resulting graph will have diameter more than $3$). Fan-structures in $(2)$ of Fig.11 or Fig.12 can be replaced by any $Fan_n$ with $n\geq 1.$ We can give the resulting graphs rainbow coloring with $4$ colors by the same method as shown in the graph $(2)$ of Fig.11 or Fig.12. Thus their rainbow connection number are at most $4.$

\begin{figure}[ht]
\begin{center}
\includegraphics[width=10cm,totalheight=35mm]{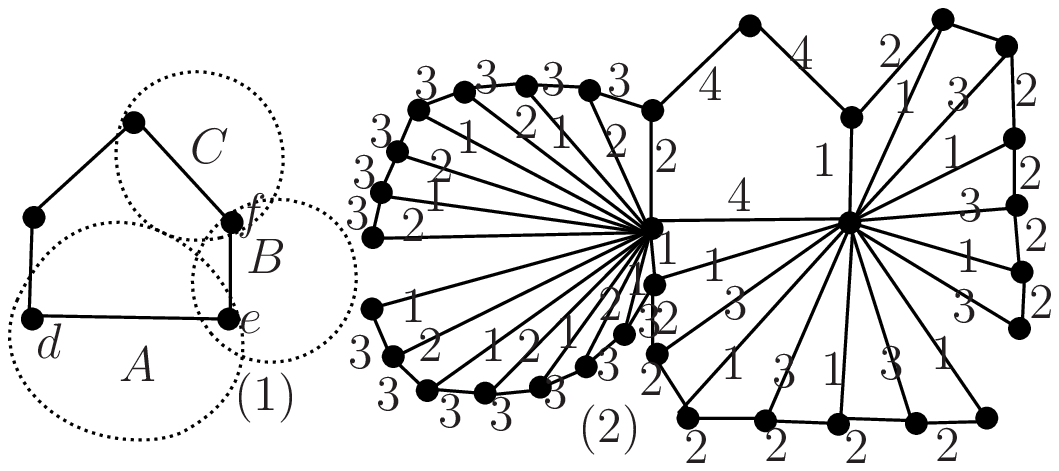}
\end{center}
\caption{$\mathrm{5-cycle ~maximal~ examlple }$}\label{fig1}
\end{figure}

%\begin{center}
%\scalebox{0.6}{\includegraphics{5-cyclemaximalexamlple.eps}}
%\end{center}

 %\begin{center}
%\scalebox{0.6}{\includegraphics{5-cycle+1Fan.eps}}
%\end{center}

\begin{figure}[ht]
\begin{center}
\includegraphics[width=10cm,totalheight=38mm]{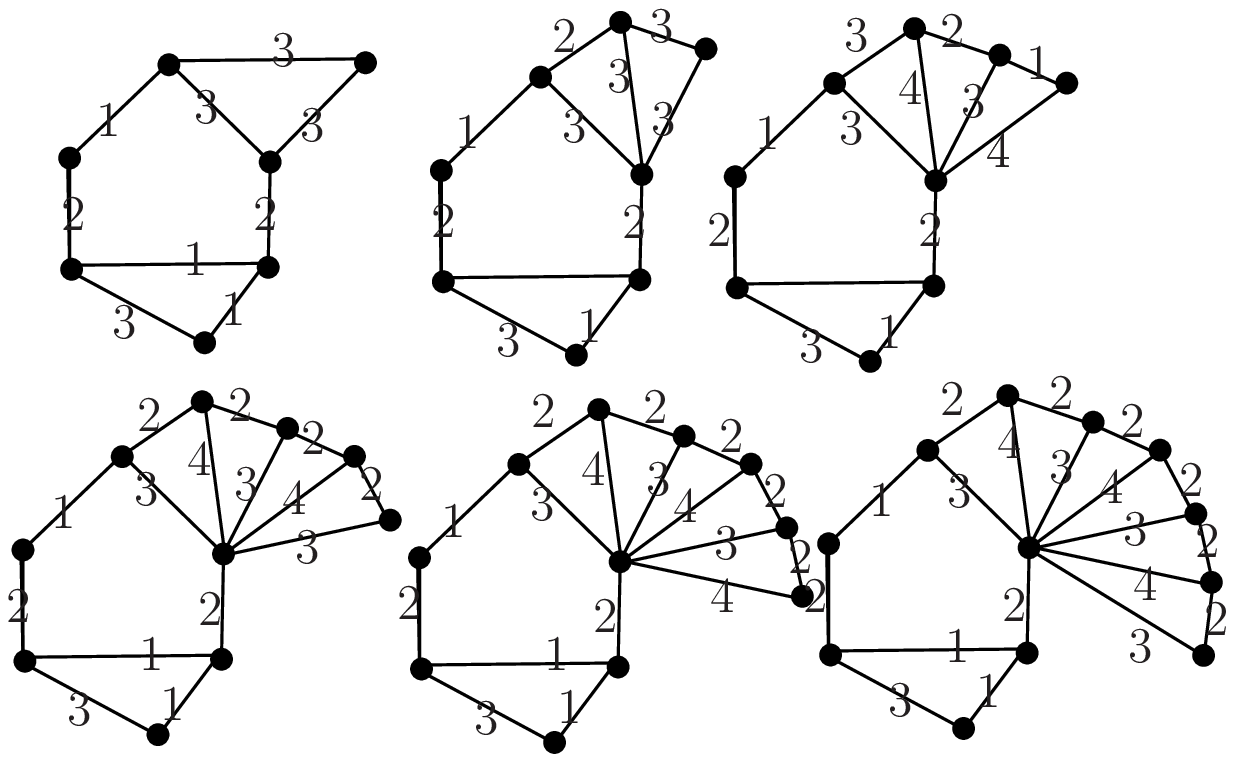}
\end{center}
\caption{$\mathrm{5-cycle+1Fan }$}\label{fig1}
\end{figure}

\noindent{\bf Case $5.$} The longest induced cycle in $G$ is $C_3,$ then it is an $MOP.$
\vspace{1mm}

Recall the construction of MOPs from Case 3 proving Theorem 2.3. When we add new vertices to Fig.13,
we can obtain a big graph having Fan-structures with $a,b$ or $c$ as their central vertices.

%\begin{center}
%\scalebox{0.6}{\includegraphics{diameter3basemaximalouterplanar1.eps}}
%\end{center}

\begin{figure}[ht]
\begin{center}
\includegraphics[width=10cm,totalheight=26mm]{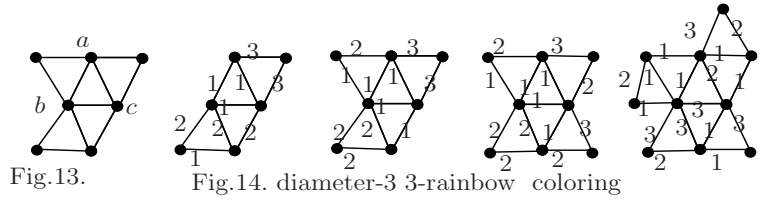}
\end{center}
\caption{$\mathrm{Diameter-3+base ~maximal~ outerplanar }$}\label{fig1}
\end{figure}

%\begin{center}
%\scalebox{0.6}{\includegraphics{diameter3+basemaximalouterplanar.eps}}
%\end{center}

%\begin{center}
%\scalebox{0.6}{\includegraphics{diameter-3+a9b9c9.eps}}
%\end{center}

\begin{figure}[ht]
\begin{center}
\includegraphics[width=8cm,totalheight=26mm]{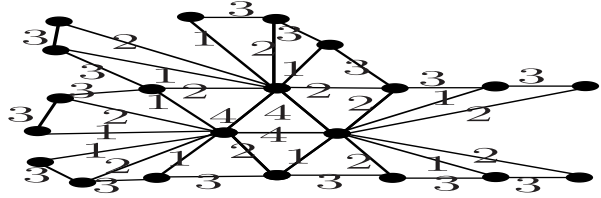}
\end{center}
\caption{$\mathrm{Diameter-3+a9b9c9}$}\label{fig1}
\end{figure}

Fig.14 gives some $MOP$s having 3-rainbow coloring and diameter $3.$ And their Fan-structures have at most $6$ vertices.
Fig.15 gives an example of MOP with diameter $3,$ where the Fan-structures have $9$ vertices. And the Fan-structures can be replaced by
any $Fan_n$ with $n\geq 1,$ and resulting graphs can be rainbow colored by $4$
colors with the same method as $Fig.15$ shown. Therefore, all graphs have rainbow connection number at most $4$ in this case.

So far we have proven the first part of the theorem. To prove the bound is sharp, it suffices to verify that
rainbow connection numbers of graph $(4)$ in Fig.16 are $4;$ which is done as follows.

We first prove that rainbow connection number of the graph $(1)$ in Fig.16 is $4.$ Consider the vertex pairs which have distance 3, and
give the only $a,b$ 3-distance path a rainbow coloring as shown in the graph $(2)$ of Fig.16. Then give the only $a,c$ 3-distance path all
possible rainbow coloring shown in $(2),(3)$ of Fig.16. In $(2),(3)$ of Fig.16, $c,e$ have the only 3-distance path $cgfe.$ Therefore,
edges $fe$ and $df$ must be colored by $1.$ Thus $d,a$ cannot be rainbow connected if we color the graph $(1)$ of Fig.16 with 3 colors.
In addition, it has a $4$ rainbow coloring as $A$ shown in $(4)$ of Fig.16. So $(1)$ of Fig.16 has rainbow connection number $4.$

We can prove that the induced subgraph $A$ of $(4)$ in Fig.16 has rainbow connection number $4 $ by the above discussion.
  There are two vertices of $A$ couldn't  be rainbow connected through the vertices of $(4)-A$ by $3$ colors which shows that $(4)$ of Fig.16 has rainbow connection number at least $4.$ However, $(4)$ of Fig.16 has a rainbow coloring by $4$ colors.
$\hfill \Box$

\begin{figure}[ht]
\begin{center}
\includegraphics[width=10cm,totalheight=45mm]{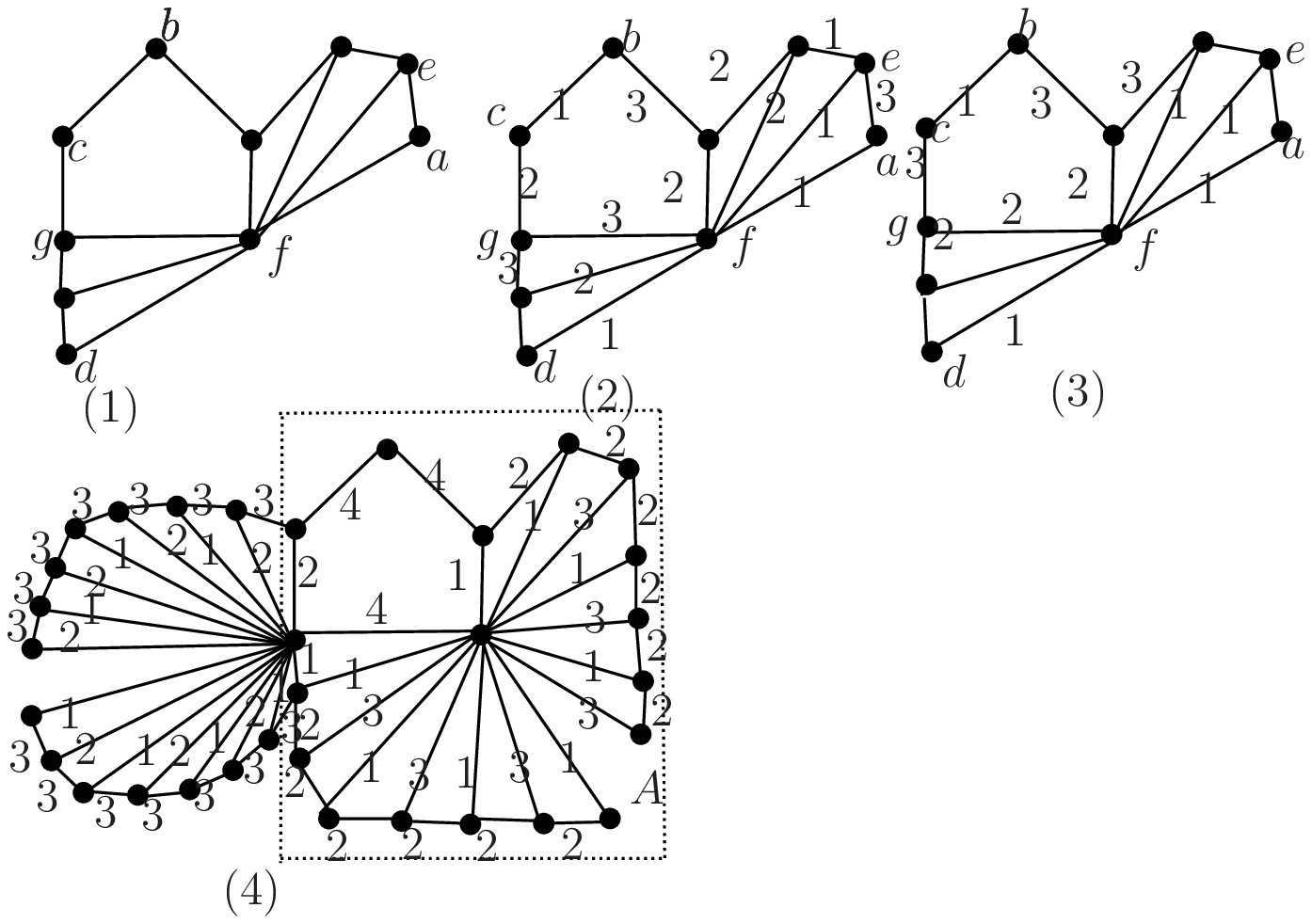}
\end{center}
\caption{$\mathrm{Optimal example}$}\label{fig1}
\end{figure}

%\begin{center}
%\scalebox{0.6}{\includegraphics{.eps}}
%\end{center}

\section{Concluding remarks}
    Uchizawa et al. \cite{UA} gave a result: Given an edge-colored graph G, decide whether G is rainbow connected, which is NP-complete even for outerplanar graphs. Thus it is interesting to study the bounds of the rainbow connection number of outerplanar graphs. In this paper, we proved that the bridgeless outerplanar graphs with diameter $2$ having rainbow connection number 2 and 3; that all the outerplanar graphs with diameter 3 have rainbow
connection number no more than 4. Moreover the bound is tight.

J. Lauri $\cite{JL}$ proved that, for outerplanar graphs given an edge-coloring, to verify it is rainbow connected is {\it NP-complete} under the coloring. Therefore it is meaning for studying approximation algorithm for MOP and outerplanar or planar graphs.

In the future, we will use the impact of $cycle$ and $Fan-structure$ to outerplanar graphs, MOP and planar graphs to design approximation algorithm for their rainbow connection numbers.

The rainbow connection comes from the communication of information between agencies of government. Our
results may help to design secure information communication networks.

\end{document}